\documentclass[12pt,leqno]{amsart}
\usepackage{amsmath,amssymb,amsfonts}
\usepackage{eucal,graphicx}

\newcommand \comment[1]{}			%  Silent version.

\newtheorem{lem}{Lemma}

\newtheorem{prop}{Proposition}
\newtheorem{thm}{Theorem}
\newcommand \mylabel[1]{\label{#1}\comment{{\rm \{#1\} }}}
\newcommand \myref[1]{\ref{#1}\comment{{\{#1\}}}}

\renewcommand{\phi}{\varphi}
\newcommand\setm{\smallsetminus}
\newcommand\inv{^{-1}}
\newcommand\textb{\text{\rm b}}

\newcommand\chib{\chi^{\textb}}
\newcommand\chiset{\chi^{\mathrm{set}}}
\newcommand\Fix{\operatorname{Fix}}
\newcommand \ind{{:}}
\newcommand\clos{\operatorname{cl}}
\newcommand\hclos{\operatorname{hcl}}
\newcommand\eset{\varnothing}

\newcommand\cH{\mathcal H}

\newcommand\cP{\mathcal P}
\newcommand \fG{\mathfrak G}
\newcommand \fH{\mathfrak H}
\newcommand \fS{\mathfrak S}
\newcommand\bbN{\mathbb N}

\allowdisplaybreaks

\begin{document}

\title{Totally frustrated states in the chromatic theory of gain graphs}
\author{Thomas Zaslavsky}
\address{Department of Mathematical Sciences\\ Binghamton University (SUNY)\\ Binghamton, NY 13902-6000 \\ U.S.A.}
\email{zaslav@math.binghamton.edu}

\begin{abstract}
We generalize proper coloring of gain graphs to totally frustrated states, where each vertex takes a value  in a set of `qualities' or `spins' that is permuted by the gain group.  (An example is the Potts model.)  
The number of totally frustrated states satisfies the usual deletion-contraction law but is matroidal only for standard coloring, where the group action is trivial or nearly regular.  
One can generalize chromatic polynomials by constructing spin sets with repeated transitive components.
\end{abstract}

\date{6 July--1 August 2006.  This version \today}
 
\maketitle
 
%%%%%%%%%%%%%%%%%%%%%%%%%%%%%%%%%%%%%%%

\section*{States, colorations, and all that}
 
When coloring a graph properly, so that the endpoints of each edge have different colors, it makes no difference what the colors are; all that matters is the number of colors in the color set.  
When coloring a gain graph, where the edges are labelled by elements of a group, that is no longer true.  The group must have a permutation action on the set of colors if the concept of a proper coloration is to mean anything, and the exact way the group acts matters very much.  
One of the major properties of graph coloring, that the number of proper colorations in $\lambda$ colors is a polynomial function of $\lambda$ that depends mainly on the graphic matroid, holds for gain graphs only when the group action on the color set is trivial, or regular, or nearly regular; this constrains $\lambda$ to take on only a fraction of all positive integer values.  
Despite this, an even more basic coloring property, the law of deletion and contraction, holds good for every color set with any group action, and there is a generalized, though non-matroidal, chromatic polynomial.

\subsection*{A state of total frustration}
A \emph{gain graph} is a graph with a function $\phi$ that assigns to each oriented edge $e$ an element $\phi(e)$ of a group called the \emph{gain group}, in such a way that reorienting the edge inverts the gain.  
A \emph{state} $s$ of the gain graph (introduced in \cite[Section 5]{CBA}) is an assignment to each vertex of an element of some set $Q$ upon which the gain group acts; $Q$ is called the set of \emph{qualities} (from \cite{RR}) or \emph{spins} (in physics).  
With a gain graph and a state, we can classify the edges as \emph{satisfied} or \emph{frustrated}: the former if, taking the edge $e$ to be oriented from vertex $v$ to vertex $w$, the equation $s_w = s_v \phi(e)$ is satisfied, and the latter if the equation is unsatisfied.  
What has been studied heretofore in connection with states has been principally the question of whether a state is satisfied (i.e., has no frustrated edges) and, if not, just how frustrated it is.  However, if we turn to states in which \emph{no} edge is satisfied, we discover a generalization of a classic problem of graph theory, the problem of proper coloring.  Our objective is to examine coloring of gain graphs from the point of view of these \emph{totally frustrated states}.

\subsection*{Properly colored}
What is new is that the set of spins is arbitrary.  In the standard theory of gain-graph coloring, from \cite[Section 4]{BG3} (the source for all properties cited herein), the color set consists of $k$ copies of the gain group $\fG$ and an extra fixed point: it is 
\[
C_k := C_k^* \cup \{0\}, \text{ where } C_k^* := \fG \times [k] 
\]
and $k$ is a nonegative integer, with $[k] := \{1,2,\ldots,k\}$ (which is void if $k=0$).  
The number of proper colorations is a polynomial function $\chi_\Phi(\lambda)$ of $\lambda = |C_k| = k|\fG|+1$, naturally called the \emph{chromatic polynomial}, that satisfies the standard deletion-contraction relation 
\begin{equation}\mylabel{E:dc}
f(\Phi) = f(\Phi\setm e) - f(\Phi/e)
\end{equation}
 for all edges $e$.  The number of proper colorations with colors taken only from $C_k^*$ is another polynomial function $\chib_\Phi(\lambda)$, the \emph{zero-free chromatic polynomial}, where now $\lambda = |C_k^*| = k|\fG|$.  
The zero-free chromatic polynomial obeys the deletion-contraction rule for edges that are not loops, and its value is not changed by the deletion of nonidentity loops.  

We want to relax the definition by admitting any spin set $Q$, and find out which properties are preserved and which are lost.

\subsection*{Strange coloring}
 The example that inspired this thought is set coloring \cite{SetCol}.  Suppose the gain group is $\fS_k$, the group of permutations of $[k]$, and the spin set $Q$ is the class of subsets of $[k]$.  
 A \emph{proper set coloration} is an assignment to each vertex of a subset $S_v$ in such a way that $S_w \neq S_v \phi(e)$ for every edge $e$, $v$ and $w$ being the endpoints of $e$.  Thus, it is a totally frustrated state of the gain graph that has an edge of every possible gain between each pair of adjacent vertices; this graph is called the \emph{$\fS_k$-expansion} of $\Delta$ and is written $\fS_k\Delta$.  
Let $\chiset(k)$ be the number of ways to do this.  This quantity is not a polynomial in any of $k$, $|Q| = 2^k$, or $|\fS_k| = k!$, so we lose something from the standard theory.  Not all is lost, however.  There is still a deletion-contraction property, so $\chiset(k)$ is what is called a Tutte invariant of gain graphs.  Our first theorem is that this is true for any group and any finite set of qualities.

\subsection*{Potts}
An example---in fact, it is an example of zero-free gain graph coloring---is the Potts model, which abstracts a partially disordered physical system such as a spin glass.  
There is a graph $\Delta$ in which each edge is marked positive or negative.  There is also a set of spins, with which we can form a state $s: V \to Q$.  
A positive edge is satisfied when it has the same spin at both ends; a negative edge is satisfied when its endpoints have different spins.  
A state has an `energy' which is a decreasing function of the number of satisfied edges.  One of the important questions is to find a lowest-energy state, or the value of the lowest energy, and especially whether there exists a completely satisfied state.  
(This is very abbreviated.  For a proper exposition with only positive edges see \cite[Section 4.4]{WelshKnots}.  The generalization to two kinds of edges is found in the physics literature and also in \cite{DM} as interpreted in \cite{BSG}.) 

To turn the Potts model into a gain graph, assume $Q$ is a group with identity element $1$.  
The Potts gain graph $\Phi$ has an edge with gain $1$ where $\Delta$ has a negative edge and it has edges with all nonidentity gains wherever $\Delta$ has a positive edge. 
A lowest-energy state of the Potts model is a state with the most frustrated edges in $\Phi$; the Potts model is satisfied when $\Phi$ is totally frustrated; and the number of frustrated edges in the Potts model is the number of satisfied edges of $\Phi$.  And in particular, the number of ways to satisfy the Potts model is the number of zero-free proper 1-colorations of $\Phi$, i.e., the value of $\chib_\Phi(|\fG|)$.

%%%%%%%%%%%%%
\section*{General theory of total frustration}

\subsection*{Technical basis}
A graph may have loops and multiple edges.  All our graphs have finite order $|V|$.  A \emph{link} is an edge that is not a loop.  
The chromatic polynomial of $\Gamma$ is written $\chi_\Gamma(\lambda)$.  
The standard closure operator on the edge set of a graph is 
\[
\clos A := A \cup \{e \notin A : \text{ the endpoints of $e$ are connected in $A$}\}.
\]
This is the closure operation of the graphic matroid (or `cycle matroid') $G(\Gamma)$.  
The notation $e\ind vw$ means that $e$ is an edge whose endpoints are  $v$ and $w$, which are equal for a loop.  If $e$ needs to be oriented (e.g.\ when evaluating its gain), the notation implies an orientation from $v$ to $w$.  

A gain graph $\Phi = (\Gamma,\phi)$ consists of an underlying graph $\Gamma = (V,E)$ and an orientable function $\phi: V \to \fG$, where $\fG$ is the gain group.  
We call $\phi$ the \emph{gain function} and $\phi(e)$ the \emph{gain} of $e$.  
By calling $\phi$ orientable we mean that its value depends on the direction of $e$ and if the direction is reversed, the gain $\phi(e)$ is inverted.  Symbolically, letting $e\inv$ denoted $e$ with the opposite orientation, $\phi(e\inv) = \phi(e)\inv$.  
The \emph{restriction} of $\Phi$ to an edge set $A$, written $\Phi|A$, is the result of deleting all edges not in $A$.  The \emph{edge-induced subgraph} $\Phi\ind A$ is $A$ together with all vertices that are endpoints of edge of $A$, and no more.  
A gain graph, or an edge set in it, is \emph{balanced} if every simple closed walk has gain, obtained by multiplying the gains of its edges in cyclic order, equal to the identity.   The number of connected components of $\Phi$ that are balanced is written $b(\Phi)$.  
  
The principal matroid in this work is the bias, or frame, matroid $G(\Phi)$ \cite{BG2}.  Its points are the edges of $\Phi$ and its rank function is $r(A) = |V| - b(\Phi|A)$.  
The class of flats determines $G(\Phi)$, of course.  
The class of balanced flats is a geometric semilattice \cite{WW} that determines what I call the \emph{balanced semimatroid} of $\Phi$, which may be defined as the class of balanced edge sets with rank as in $G(\Phi)$.  
In general one may think of a semimatroid $S$ as the class of sets in a matroid $M$ on $E(S) \cup \{e_0\}$ whose closures do not contain $e_0$, together with the rank function on these sets.  We call $M$ the \emph{completion} of the semimatroid.  According to \cite[Theorem 3.2]{WW}, the completion is unique.  
The completion of the balanced semimatroid of $\Phi$ is the complete lift matroid $L_0(\Phi)$ \cite{BG2}, whose contraction $L_0(\Phi)/e_o$ equals $G(\Gamma)$.

A \emph{spin set} is a set $Q$ upon which there is a right action of $\fG$.  The action is \emph{trivial} if every spin is fixed by every group element and \emph{regular} if only the identity element has any fixed points.  
A \emph{state} of $\Phi$ is any function $s: V \to Q$.  It is \emph{totally frustrated} if every edge is frustrated.  

We must define deletion and contraction of an edge.  Deletion is obvious.  To contract a link we need \emph{switching}.  A \emph{switching function} $\eta: V \to \fG$ gives a switched graph $\Phi^\eta$ whose underlying graph is the same as that of $\Phi$ and whose gain function is $\phi^\eta$, defined by $\phi^\eta(e) := \eta_v\inv \phi(e) \eta_w$ for any edge $e\ind vw$.  It is always possible to switch so a given link has gain $1$, the group identity, and indeed so that the gains on a chosen forest are all $1$.   In order to contract a link $e$ we switch so it does have gain $1$, then we delete it and identify its endpoints.  The gains do not change except in the switching step.  
Contraction of a loop will not be needed; for it see \cite{BG1}.  
Note that $\Phi/e$ is uniquely defined only up to switching.

There is one more aspect of switching that is essential: switching acts on states as well as gains.  We define $s^\eta$ by $s^\eta_v := s_v\eta_v$; in words, a switching function acts on a state in the obvious way.  The point is that the set of satisfied edges, $I_\Phi(s)$, remains the same: 
\[
I_{\Phi^\eta}(s^\eta) = I_\Phi(s).
\]
%CALCULATION IN CASE WANTED BY REFEREE:  \[ s_w = s_v \phi(e) \iff s_w \eta_w = (s_v \eta_v) \phi^\eta(e) \iff s^eta_w = s_v^eta \phi^\eta(e). \]
Obviously, therefore, the number of totally frustrated states is unaffected by switching and we may assume that $\Phi$ is already switched so that $\phi(e) = 1$.
%(Switching is wonderful!  Credit to J.J.\ Seidel, who introduced it, implicitly, for signed graphs.)

\subsection*{States vs.\ colorations}
The difference between a state, with an arbitrary spin set, and a coloration, whose spin set (or `color set') is $C_k$ or $C_k^*$, is that in a coloration the spin set yields properties very similar to those of ordinary graph coloring.  For instance, the set of frustrated edges in a coloration is closed in the frame matroid $G(\Phi)$ \cite{BG3}.  

One could say that the crux of the difference is the behavior of loops---not surprisingly, because knowing which spins on its supporting vertex satisfy a loop with gain $g$ is the same as knowing the fixed points of $g$ acting on $Q$, and by Theorems \myref{T:chromatic} and \myref{T:nonmatroid}, that is what decides whether or not $\chi_\Phi(Q)$ does or does not equal the chromatic or zero-free chromatic polynomial of $\Phi$ (or the chromatic polynomial of the underlying graph).  
The proofs of those theorems show that the most basic question about loops is whether the number of totally frustrated states of a nonidentity loop is affected by the exact gain, or in other words, whether every $g\neq1$ fixes the same number of spins.

\subsection*{Equivalence-class colorations}
This is an example in which we have a finite spin set $Q$ partitioned into $Q_1,\ldots,Q_r$.  The gain group is the group of permutations of $Q$ that respect the partition, that is, $\fG$ is the product of the symmetric groups of the subsets $Q_i$.  Suppose we have a graph $\Delta$.  If we take $\Phi = \fG\Delta$, in which each edge of $\Delta$ is replaced by edges having every gain in $\fG$, then a totally frustrated state of $\Phi$ is precisely a state in which adjacent vertices of $\Delta$ have inequivalent spins.

\subsection*{The chromatic function}
Choosing a spin set $Q$, the \emph{$Q$-chromatic function} of $\Phi$ is
\[
\chi_\Phi(Q) := \text{ the number of totally frustrated states}.
\]
This is a finite number if $Q$ is finite.

\begin{prop}\mylabel{L:bal}
If $\Phi$ is balanced, then $\chi_\Phi(Q) = \chi_\Phi(|Q|) = \chi_\Gamma(|Q|)$.
\end{prop}

\begin{proof}
By switching we may assume all gains equal $1$.  Clearly, then $\chi_\Phi(Q) = \chi_\Gamma(|Q|) = \chi_\Phi(|Q|)$.
\end{proof}

\begin{thm}\mylabel{T:dc}
If $Q$ is finite, the $Q$-chromatic function of $\fG$-gain graphs of finite order has the deletion-contraction property \eqref{E:dc} with respect to all links $e$.
\end{thm}

\begin{proof}
We simply classify the totally frustrated states of $\Phi\setm e$ according to whether $e$ is or is not frustrated.  A state for which $e$ is frustrated is a totally frustrated state of $\Phi$.  The criterion for $e$ to be satisfied is that its endpoints have the same spin.  Hence a state in which $e$ is satisfied contracts to a totally frustrated state of $\Phi/e$, and conversely, any totally frustrated state of $\Phi/e$ defines a unique state of $\Phi$ in which $e$ and only $e$ is satisfied.  This proves the theorem.
\end{proof}

\subsection*{Decomposition}
A normalization of the $Q$-chromatic function is $\theta(\Phi) := |Q|^{-b(\Phi)} \chi_\Phi(Q)$.  The same normalization applied to the chromatic polynomial, i.e., $\lambda^{-b(\Phi)} \chi_\Phi(\lambda)$, gives the characteristic polynomial of $G(\Phi)$ \cite[Section 5]{BG3}.

\begin{prop}\mylabel{T:mult}
Assume finite $Q$ and $\Phi$ and suppose $\Phi'$ and $\Phi''$ are subgraphs whose union is $\Phi$.  
If they are disjoint, or if their intersection is a single vertex and at least one of them is balanced, then 
\[ 
\theta(\Phi)= \theta(\Phi') \theta(\Phi'').
\]
\end{prop}

\begin{proof}
If $\Phi'$ and $\Phi''$ are vertex disjoint, then multiplicativity is obvious.  From this one can see that it suffices to assume $\Phi'$ and $\Phi''$ are connected.

Suppose the intersection is a vertex $v$ and $\Phi'$ is balanced.  
A state of $\Phi$ is totally frustrated if and only if it is assembled from a totally frustrated state $s''$ of $\Phi''$ and a totally frustrated state $s'$ of $\Phi'$ that agrees with $s''$ on $v$.  The question is how the number of such $s'$ depends on $s''_v$.  The number is independent of $s''_v$, indeed it equals $\chi_{\Phi'}(|Q|)/|Q|$ (by switching as at Proposition \myref{L:bal}), hence it equals $\theta_{\Phi'}(|Q|)$.  Multiplicativity follows.
\end{proof}

A particular case is the obvious fact that the chromatic function is multiplicative on connected components: if $\Phi$ has components $\Phi_1, \Phi_2, \ldots, \Phi_m$ then 
$$
\chi_\Phi(Q) = \chi_{\Phi_1}(Q) \chi_{\Phi_2}(Q) \cdots \chi_{\Phi_m}(Q) .
$$
This, together with the deletion-contraction law and the facts that $\chi_\Phi(Q)$ is an isomorphism invariant and $\chi_\eset(Q) = 1$, means that the $|Q|$-chromatic function satisfies the definition of a Tutte invariant of gain graphs, thus being another in a long list of such invariants.  (We forgo any further discussion of Tutte invariance here; see \cite{BOTutte}.)

\subsection*{Chromatic or not chromatic}
To say that $\chi_\Phi(Q)$ is an \emph{evaluation} of a function $\chi(\lambda)$ means there is a fixed value $\lambda_0$ such that $\chi_\Phi(Q) = \chi(\lambda_0)$ for every $\fG$-gain graph $\Phi$.

\begin{thm}\mylabel{T:chromatic}
Assume $Q$ is a finite set of spins.  

{\rm(a)}  If there is $q_0 \in Q$ such that every nonidentity group element has fixed set $\{q_0\}$, then $\chi_\Phi(Q)$ is the evaluation of the chromatic polynomial $\chi_\Phi(\lambda)$ at $\lambda = |Q|$.  If the assumption fails, then $\chi_\Phi(Q)$ is not an evaluation of $\chi_\Phi(\lambda)$.

{\rm(b)}  If $\fG$ acts regularly on $Q$, then $\chi_\Phi(Q)$ is the evaluation of the zero-free chromatic polynomial $\chib_\Phi(\lambda)$ at $\lambda = |Q|$.  If the assumption fails, then $\chi_\Phi(Q)$ is not an evaluation of $\chib_\Phi(\lambda)$.

{\rm(c)}  If the action of $\fG$ is trivial, then $\chi_\Phi(Q) = \chi_\Gamma(|Q|)$, the evaluation at $|Q|$ of the chromatic polynomial of the underlying graph.  If the action is nontrivial, then $\chi_\Phi(Q)$ is not an evaluation of $\chi_\Gamma(\lambda)$.
\end{thm}

\begin{proof}
We prove the implication in part (a) in stages.  The underpinning is that the chromatic polynomial satisfies deletion-contraction for all links.  Thus, if we prove the theorem for graphs without links, it follows by induction on the number of edges using Theorem \myref{T:dc}.  
The chromatic polynomial and the $Q$-chromatic function both equal zero when $\Phi$ has an identity loop, so we may assume $\Phi$ has no edges other than nonidentity loops.
Furthermore, both chromatic polynomial and $Q$-chromatic function are multiplicative on connected components, so we may assume $\Phi$ is connected.  That is, $\Phi$ has a single vertex with some number of nonidentity loops.

If there are no loops, both $\chi_\Phi(Q)$ and $\chi_\Phi(\lambda)$ equal $1$.  If there is at least one loop, then $\chi_\Phi(\lambda) = \lambda-1$.  Now, let $G$ be the set of gains of the loops of $\Phi$.  To be totally frustrated, a state $s$ must have $s_v \notin \Fix(g)$, the fixed set of $g$, for every $g \in G$.  The only way that can give $\lambda-1$ is if $\Fix(g)$ is the same set $F$ for every nonidentity element of the gain group and $\lambda = |Q| - |F|+1$.

We can evaluate $\lambda$, thereby determining $|F|$.  Consider $\Phi$ with vertex set $\{v,w\}$ and one edge that is a link $e\ind vw$.  Then $\chi_\Phi(Q) = |Q|(|Q| - 1)$.  This should be the chromatic polynomial evaluated at $\lambda$, but the chromatic polynomial is $\lambda(\lambda-1)$.  Thus, $\lambda = |Q|$.  It follows that $|F|=1$.  

So, we have necessary conditions for the $Q$-chromatic function to be an evaluation of the chromatic polynomial, but the proof also shows their sufficiency.  That concludes the proof of part (a).

The proof of part (b) is similar.  A nonidentity loop is never satisfied so it can be discarded without altering the number of totally frustrated states.

For part (c), note that if the action is trivial, then the gains do not matter.  Conversely, if there is a $g$ with nontrivial action, consider the gain graph with one vertex and one loop, whose gain is $g$.  Then $\chi_\Gamma(\lambda) = 0$ but $\chi_\Phi(Q) = |Q|$.
\end{proof}

Theorem \myref{T:chromatic} demonstrates that the $Q$-chromatic function equals the chromatic polynomial only when $Q$ is essentially a $C_k$.  
Take $Q$ as in Theorem \myref{T:chromatic}(a), delete $q_0$, and divide up the rest into orbits of $\fG$.  Each orbit can be identified with $\fG$ (with the right regular action) since there are no fixed points of any group element other than the identity.  Thus $|Q|$ is $C_k$ in disguise, $k$ being the number of orbits in $Q \setm q_0$.
Similar remarks hold good for the second part of the theorem, with color set $C_k^*$.
In fact the $Q$-chromatic function is not even determined by the matroid, or semimatroid, unless $Q$ has exactly the form stated in Theorem \myref{T:chromatic}; this will be Theorem \myref{T:nonmatroid}.

\subsection*{Matroid invariance}
We can strengthen the second halves in Theorem \myref{T:chromatic}(a, b, c) to a characterization of when the number of totally frustrated states is a matroid or semimatroid invariant.  
The matroid involved is the frame matroid $G(\Phi)$.  

\begin{thm}\mylabel{T:nonmatroid}
Let $Q$ be a finite spin set.

{\rm(a)}  The $Q$-chromatic function, as a function of the gain graph, is determined by the frame matroid $G(\Phi)$ and the numbers of components and balanced components of $\Phi$, if and only if $Q$ meets the conditions of Theorem \myref{T:chromatic}(a).

{\rm(b)}  The $Q$-chromatic function, as a function of the gain graph, is determined by the balanced semimatroid of $\Phi$ and the numbers of components and balanced components of $\Phi$ if and only if $\fG$ acts regularly or trivially upon $Q$, as in Theorem \myref{T:chromatic}(b) or (c).
\end{thm}

\begin{proof}[Proof of (a)]
If $Q$ does meet the conditions, then it is an evaluation of $\chi_\Phi(\lambda)$ (by Theorem \myref{T:chromatic}(a)), which in turn is equal to $\lambda^{b(\Phi)}$ times the characteristic polynomial of $G(\Phi)$ \cite[Section 5]{BG3}.

For the converse, suppose the $Q$-chromatic function is determined by the stated information.  Let $g,h\in\fG$, both not the identity.  

The gain graph $\Phi_g$ that consists of one vertex and one loop with gain $g$ has matroid isomorphic to a coloop.  The chromatic function is $|Q| - |\Fix(g)|$.  Since $\Phi_h$ has the same matroid and component numbers as $\Phi_g$, it must have the same chromatic function.  It follows that every group element other than $1$ must have the same number $f$ of fixed points.

The gain graph $\Phi_{g,h}$ has vertex $v_1$ with a loop of gain $g$ and $v_2$ with a loop of gain $h$ and a link $e\ind v_1v_2$ with gain $1$.  At $v_1$ the spin is $q_1 \notin \Fix(g)$, and at $v_2$ the spin is $q_2 \notin \{q_1\} \cup \Fix(h)$.  The chromatic function is 
\begin{align*}%\mylabel{E:}
\chi(Q) &= \sum_{q_1 \in Q\setm\Fix(g)} \big[ |Q| - |\{q_1\} \cup \Fix(h)| \big] \\
	&= \sum_{q_1 \notin \Fix(g)\cup\Fix(h)} \big[ |Q| - f - 1 \big] + \sum_{q_1 \in \Fix(h)\setm\Fix(g)} \big[ |Q| - f \big] \\
	&= (|Q|-f)(|Q|-f-1) + |\Fix(h)\setm\Fix(g)| .
\end{align*}
This value cannot depend on the choices of $g,h\neq1$ because those do not change the matroid.  Since taking $g=h$ gives value $(|Q|-f)(|Q|-f-1)$, it follows that every nonidentity element has the same fixed set.

In effect, $Q$ is the disjoint union $\fG\times[k_1] \cup Q_2\times[k_2]$ where $|Q_2| = 1$.  (That is because the nontrivial orbits of $Q$ have no fixed points of any nonidentity element of $\fG$.)  
%This is the situation of Theorem \myref{T:poly} with $p=2$, $Q_1=\fG$, and $|Q_2|=1$.  
Let $\fG$ act on $Q$ by its effect on the first component, i.e., the element of $Q_i$.  We compare two gain graphs.

Consider first the gain graph $\Phi_2$ that has vertices $v_1$ and $v_2$ with two links joining them, one having gain $1$ and the other with gain $g\neq1$.  To get a totally frustrated state we choose spin $q_1$ for $v_1$.  There are $|Q|-1$ choices for $q_2$ at $v_2$ if $q_1 \in Q_2\times[k_2]$ but $|Q|-2$ choices if $q_1 \in Q_1\times[k_1]$, since that $q_1$ implies $q_2 \neq q_1$ and $q_1g$.  
Thus, $\chi_{\Phi_2;Q_1,Q_2}(k_1,k_2) = (|Q|-1)^2 + k_2-1$.  

Second, consider $\Phi_1$ that has the same vertices and a link with gain $1$, but also a loop at $v_1$ with gain $g\neq1$.  The matroid and the numbers of components and balanced components are the same, but it is easy to see that $\chi_{\Phi_1;Q_1,Q_2}(k_1,k_2) = (|Q|-2)^2$.  We conclude that the only case in which the $|Q|$-chromatic function can be determined by the stated information is that in which $k_2=1$.
\end{proof}

\begin{proof}[Proof of (b)]
If $Q$ does meet the conditions, then it is an evaluation of $\chib_\Phi(\lambda)$ (by Theorem \myref{T:chromatic}(b,c)), which in turn is equal to $\lambda^{b(\Phi)}$ times the characteristic polynomial of the semilattice of balanced flats \cite[Section 5]{BG3}.

The proof of the converse proceeds in the same way as in part (a) to establish the form of $Q$, i.e., $p=2$, $Q_1=\fG$, and $|Q_2|=1$, since the data of part (b) agree for all the graphs we compared in the steps leading to that conclusion.

The final step is much more complicated than in part (a) because we need nonisomorphic gain graphs that have the same balanced sets and their ranks.  That is impossible with only two vertices.  Figures \myref{F:balanceisoma} and \myref{F:balanceisomb} show two gain graphs with this property.  The balanced sets are all those that do not contain a digon and, if they contain exactly one $f_i$, do not complete a circuit of four edges; these sets are the same in both graphs.

\begin{figure}\mylabel{F:balanceisoma}
\includegraphics{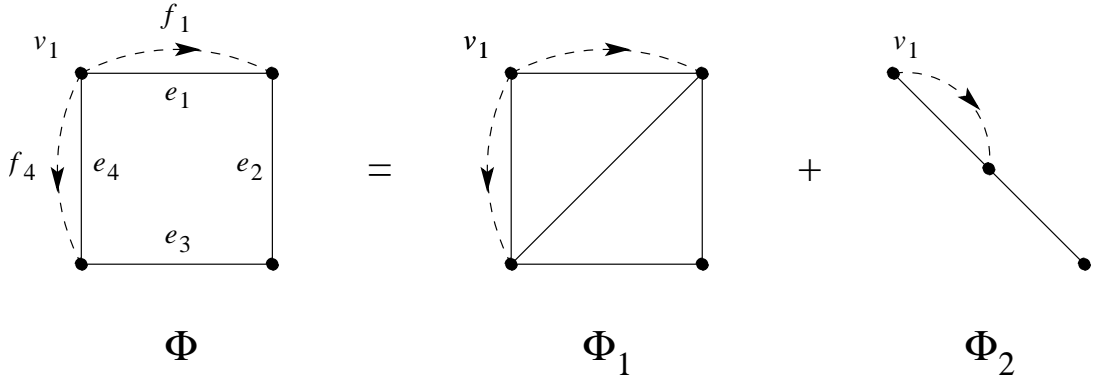}
\caption{A gain graph $\Phi$ with four identity edges $e_i$ and two adjacent edges that have gain $\phi(f_i) = g\neq1$, showing how its chromatic function decomposes by addition and contraction of an edge having gain $1$.  (In the contracted graph, multiple edges with identical gain are suppressed.)}
\end{figure}

\begin{figure}\mylabel{F:balanceisomb}
\includegraphics{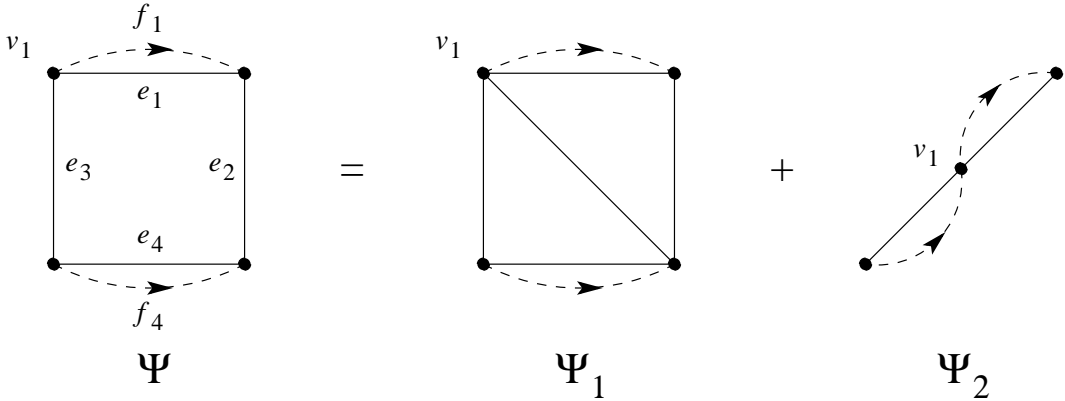}
\caption{A gain graph $\Psi$ with four identity edges $e_i$ and two nonadjacent edges that have gain $\phi(f_i) = g\neq1$.}
\end{figure}

The calculation of the chromatic function is not so simple.  I used deletion-contraction in reverse, by means of which the chromatic function of each graph is expressed as the sum of two other chromatic functions that are easier to work with (see the figures).  
I calculated the number of totally frustrated states of each of the four graphs by the usual hand method of building the state from vertex to vertex, starting with spin $q_1$ at $v_1$ and treating separately spin $q_1 \in Q_1\times[k_1]$ and $q_1 \in Q_2\times[k_2]$.  I checked the result by comparing it, with $k_2=0$, to the zero-free chromatic polynomial computed from the semilattice of balanced flats as in \cite[Section 5]{BG3}, which by Theorem \myref{T:chromatic}(b) ought to be the same (and was).  I omit the lengthy details.  
The conclusion is that 
\begin{align*}%\mylabel{E:}
\chi_{\Phi;Q_1,Q_2}(k_1,k_2) &= \lambda(\lambda-2)[\lambda^2-4\lambda+5] + k_2 \big[2\lambda^2 - 7\lambda + 7 \big] ,      \\
\chi_{\Psi;Q_1,Q_2}(k_1,k_2) &=  \lambda(\lambda-2)[\lambda^2-4\lambda+5] + k_2 \big[2\lambda^2 - 8\lambda + 7 + k_2 \big] ,   
\end{align*}
where $\lambda = |Q|$.  
The difference between these is $k_2(\lambda-k_2)$.  Therefore, they are equal only when $k_2=0$, so the number of totally frustrated states is $\chib_\Phi(|Q|)$, or $k_1=0$, so the action of $\fG$ is trivial.  In the latter case, we know $\chi_\Phi(Q) = \chi_\Gamma(k_2)$ by Lemma \myref{L:bal}, but we have yet to prove that $\chi_\Gamma(\lambda)$ is determined by the balanced semimatroid.  This follows from the fact that the completion of the balanced semimatroid is $L_0(\Phi)$, whose contraction by the extra point $e_0$ is the graphic matroid $G(\Gamma)$.  Since the completion is unique, the balanced semimatroid determines $G(\Gamma)$; this in turn determines $\chi_\Gamma(\lambda)$ as $\lambda^{c(\Gamma)}$ times the characteristic polynomial of $G(\Gamma)$.  Hence, the semimatroid does determine the number of totally frustrated states.
\end{proof}

\subsection*{Multiplicativity not matroidal}
Proposition \myref{T:mult} is one case of a decomposition into separators of the frame matroid.  Under the assumptions, the edge sets of both subgraphs are complementary separators of $G(\Phi)$.  
It is natural to ask whether, if $\Phi$ has any subgraphs $\Phi'$ and $\Phi''$ whose edge sets are complementary matroid separators, then $\theta$ is multiplicative on them.  The answer is no.  
Suppose $\Phi$ is unbalanced and $e$ is a link such that $\Phi'' := \Phi \setm e$ is balanced; then $\{e\}$ and its complement are separators of the matroid.
Consider the specific example of $K_{|Q|+1}$ with gains all $1$ except for one edge $e$, and $Q$ with a regular action.  It is easy to compute the functions to see that $\theta(\Phi) \neq \theta(\Phi\ind e) \theta(\Phi\setm e)$.

\subsection*{Holonomy}
To get a more precise formula for the chromatic function we need a new concept.  
First we define the gain of a walk $W$: it has gain $\phi(W)$ equal to the product of the gains of its edges, in the order traversed by $W$.  For a tree $T$ and an edge $e$ whose endpoints are connected by $T$, let $T(e)$ be the component of $T$ that contains the endpoints of $e$.  
Given a vertex $v_0$ of $T(e)$, let $W_e$ be a minimal closed walk in $T_0 \cup \{e_0\}$ that starts at $v_0$ and contains $e$.

Now, take an edge set $A\subseteq E$ and let $T$ be a maximal forest in $A$.  Let $A_0$ be the component of $A$ that contains $v_0$.  For each $e \in A_0 \setm T$, let $W_e$ be a minimal closed walk from $v_0$ in $T_0 \cup \{e_0\}$ that contains $e$.  
The \emph{holonomy group} of $A$ with respect to $T$ at a base vertex $v_0$ is
\[
\fH_A(v_0,T) := \langle \phi(W_e): e \in A_0 \setm T_0 \rangle ,
\]
where the angle brackets indicate the subgroup of $\fG$ generated by the gains $\phi(W_e)$.  The generator $\phi(W_e)$ is called the \emph{holonomy of $e$}; it is $1$ if and only if the edge set of $W_e$ is balanced, as when $e\in T$.  Should it happen that $T$ has identity gain on all edges, the definition simplifies to $\fH_A(v_0,T) := \langle \phi(e): e \in A_0 \setm T_0 \rangle.$  Should $v_0$ happen to be isolated in $V$, then its holonomy group is the subgroup of $\fG$ generated by the gains of loops at $v_0$.

One can define the holonomy group in terms of contraction.  If we choose $T$ and contract it, then $\fH_A(v_0,T)$ is the subgroup of $\fG$ generated by the gains of the loops of the contracted graph that are incident with $v_0$ and belong to $A$.

\begin{lem}\mylabel{L:holotree}
The holonomy groups of $A$ with respect to different maximal forests are the same.
\end{lem}

\begin{proof}
By Lemma \myref{L:holobase} we may assume the basepoint $v_0$ is fixed and only $T$ changes to $T'$.  We show that the holonomy of $e$ with respect to one forest is contained in the holonomy group with respect to the other.  

Consider an edge $e\ind vw$ in $A_0 \setm T$.  Let $T_v$ denote the path in $T$ from $v_0$ to $v$.  Then $\phi(W_e) = \phi(T_v) \phi(e) \phi(T_w)\inv$, and similarly for $W'_e$.  Let $T'_v = v_0e_1v_1 \cdots e_lv_l$ and $T'_w = v_0f_1w_1 \cdots f_mw_m$.  Then 
\begin{align*}
\phi(W'_e) &= \phi(T'_v) \phi(e) \phi(T'_w)\inv \\
&= \phi(W_{e_1} \cdots W_{e_l} \cdot W_e \cdot W_{f_m}\inv \cdots W_{f_1}\inv ) \\
&= \prod_1^l \phi(W_{e_i}) \cdot \phi(W_e) \cdot \prod_m^1 \phi(W_{f_j})\inv \in \fH_A(v_0,T) ,
\end{align*}
so $\fH_A(v_0,T') \subseteq \fH_A(v_0,T)$; and similarly there is the opposite inclusion.
\end{proof}

Although it is no longer necessary to specify $T$ in the notation for the holonomy group $\fH_A(v)$, we shall still do so when it seems especially helpful. 

\begin{lem}\mylabel{L:holobase}
Let $A \subseteq E$ and $T$ a maximal forest in $A$.  For any two vertices in the same component of $A$, their holonomy groups with respect to $A$ and $T$ are conjugate in $\fG$.
\end{lem}

\begin{proof}
If $v$ and $w$ are the vertices, let $P$ be the path in $T$ from $v$ to $w$.  
Then the walk $W_e(w)$ based at $w$ is the reduced form of the walk $P\inv W_e(v) P$, so $\phi(W_e(w)) = \phi(P)\inv \phi(W_e(v)) \phi(P)$.  
It follows that $\fH_A(w,T) = \phi(P)\inv \fH_A(v,T) \phi(P)$.
\end{proof}

\begin{lem}\mylabel{L:holoswitch}
Fixing $A$, $T$, and $v_0$, switching $\Phi$ conjugates the holonomy group.
\end{lem}

\begin{proof}
Observe that $\phi^\eta(W_e) = \eta_v\inv \phi(W_e) \eta_v$.
\end{proof}

The consequence is that, given $A$, although the fixed set of the holonomy group may depend on the basepoint and switching, the size of the fixed set is independent of these choices as long as the basepoint stays in the same component.  That is because conjugating a subgroup $\fH$ by $\alpha \in \fG$ changes $\Fix \fH$ to $\Fix (\fH^\alpha) = (\Fix \fH) \alpha$.  Thus, we are justified in defining
\[
h_Q(A) := |\Fix(\fH_A(v))| 
\]
for a connected edge set $A$; we assume $v$ is chosen in the vertex set of $A$; and if $A$ is empty the holonomy group is trivial so $h_Q(\eset) = 1$; and $h_Q(A)$ is invariant under switching of $\Phi$.

Another consequence is the case of a balanced edge set, obtained by applying the following lemma to a connected edge set in a larger gain graph $\Phi$.

\begin{lem}\mylabel{L:balanced}
If a gain graph is connected, its holonomy group is trivial if and only if it is balanced.
\end{lem}

\begin{proof}
By switching assume the graph contains an identity spanning tree.  There is a nonidentity edge if and only if the graph is unbalanced.
\end{proof}

Next we define \emph{holonomy closure}.  Again, $T$ is some maximal forest in $A$.  The holonomy closure of $A$ is
\begin{align*}
\hclos A := \{ e \in \clos A : \phi(W_e) \in \fH_A(v,T) \} 
\end{align*}
where $v$ is some fixed vertex in $V(T(e))$.  Note that $\hclos A \supseteq A$.  
(It might help to think of the holonomy closure as the inverse image of the holonomy group $\fH_A(v)$ under the $v,T$-gain function $\phi_{v,T}(e) := \phi(W_e)$.) 
We know from Lemmas \myref{L:holotree} and \myref{L:holobase} that this closure is independent of the choices.  We mention $T$ only because $W_e$ depends on which $T$ we pick.  
The definition simplifies if $T$ happens to have all identity gains; then 
\begin{align*}
\hclos A = \{ e \in \clos A : \phi(e) \in \fH_A(v,T) \} .
\end{align*}

A set that is its own holonomy closure is, of course, called \emph{holonomy closed}.  We write $\cH_\Phi$ for the class of holonomy-closed edge sets.

\subsection*{Satisfied edge sets}
An arbitrary state $s$ has a set $I(s)$ of satisfied edges; we ask what kind of set this can be.  
We want a characterization in terms of the gains and gain group, independent of the particular actions.  The detailed formula we want for the chromatic function comes from M\"obius inversion over the sets $I(s)$; knowing what they may be tells us the poset over which to invert.

Recall that the satisfied edges are invariant under switching.

\begin{lem}\mylabel{L:frustholo}
The satisfied edges of a state constitute a holonomy-closed set.
\end{lem}

\begin{proof}
Take a state $s$ and an edge $e\ind vw$ in the holonomy closure of $I(s)$.  Choose a spanning tree $T_0$ of the component of $I(s)$ that contains the endpoints of $e$ and switch by $\eta_v := \phi((T_0)_v)\inv$ so $T_0$ has identity gain.  Then $s$ is constant on $V(T_0)$; say $s_v = q \in Q$ for every $v \in V(T_0)$.  We want to prove that the switched holonomy $\phi^\eta(W_e) = \phi(e)$ belongs to the stabilizer of $q$, $\fG_q$.  This will imply that it is in $I(s^\eta)$, which we know equals $I(s)$.  

Each holonomy generator $\phi^\eta(W_f) = \phi(f)$ lies in $\fG_q$ because $f \in I(s)$.  Therefore, $\fH^\eta_A(v) \leq \fG_q$, and that is what we need.
\end{proof}

\subsection*{A very detailed formula}
We are now ready to employ the standard method of M\"obius inversion \cite{FCT,EC1} to get an exact formula for the $Q$-chromatic function.  In ordinary coloring theory the formula is quite simple because the number of all colorations, proper or not, is simply a power of the number of colors (see \cite[Section 9]{FCT}), but in state coloring the result has to be expressed in terms of fixed sets of holonomy groups.  We state two versions of the formula.  The first has fewer terms but involves the M\"obius function of the lattice $\cH_\Phi$ of holonomy-closed sets, about which nothing is known.  The second, which is just an inclusion-exclusion formula, is simpler but has more terms.

\begin{thm}\mylabel{T:detailed}
For a finite gain graph $\Phi$ with a finite spin set $Q$,
\begin{align*}%\mylabel{E:}
\chi_\Phi(Q) &= \sum_{A \in \cH_\Phi}  \mu_{\cH_\Phi}(\eset,A)  \prod_{j=1}^m h_Q(A_j) \\
	&= \sum_{A\subseteq E}  (-1)^{|A|}  \prod_{j=1}^m h_Q(A_j) ,
\end{align*}
where $A_1,\ldots,A_m$ are the connected components of $A$.
\end{thm}

An isolated vertex must be treated as a connected component, as its value of $h_Q$ is $|Q|$.

\begin{proof}
We shall prove the first formula, summing over the class $\cH_\Phi$, but that of the second is identical except for replacing $\cH_\Phi$ by $\cP(E)$ with its M\"obius function $\mu(\eset,A) = (-1)^{|A|}$.

Let $f(A)$ be the number of states for which $I(s) \supseteq A$ and let $g(A)$ be the number such that $I(s) = A$.  Since every possible set of satisfied edges belongs to $\cH_\Phi$, we see that
\[
f(B) = \sum_{A \in \cH_\Phi: A \supseteq B} g(A)
\]
for every $B \in \cH_\Phi$.  By M\"obius inversion, 
\[
g(B) = \sum_{A \in \cH_\Phi} \mu_{\cH_\Phi}(B,A) f(A).
\]
Setting $B = \eset$ we get
\[
\chi_\Phi(Q) = g(\eset) = \sum_{A \in \cH_\Phi} \mu_{\cH_\Phi}(\eset,A) f(A).
\]

To finish the proof we have to interpret $f(A)$.  Let $T$ be a maximal forest in $A$ and switch so $T$ has identity gains.  Then any state counted by $f(A)$ is constant on each component $A_j$, having let us say spin $q_j$.  For each other edge $e \in A_j$ we must have $q_j \phi(e) = q_i$; thus, $q_j$ can be any spin that is fixed by every gain $\phi(e)$ for $e \in A_j$.  These gains are the generators of $\fH_{A_j}(v,T)$, so the possible spins $q_j$ are precisely those that lie in $\Fix \fH_{A_j}(v,T)$.  The number of these is $h_Q(A_j)$.  The value of $f(A)$ is the number of ways to choose one spin for each component, i.e., the product of all $h_Q(A_j)$.  That proves the formula.
\end{proof}

%%%%%%%%%
\section*{A grand polynomial}

Despite all the difficulties about matroids, there is a way to make the $Q$-chromatic function into a polynomial that generalizes the chromatic polynomial.  
Let us have spin sets $Q_1,Q_2, \ldots, Q_p$, that is, each is a set with a $\fG$-action, and to avoid notational difficulty suppose that all the sets $Q_i$ and $Q_i \times \bbN$ are pairwise disjoint.  ($\bbN$ is the set of nonnegative integers.)  
Write $\Fix_i \fH$ for the fixed set of the action of $\fH$ on $Q_i$, where $\fH$ is a subgroup of $\fG$.  
Set 
\[
Q := Q_{k_1,k_2,\ldots,k_p} := Q_1\times[k_1] \cup Q_2\times[k_2] \cup \cdots \cup Q_p\times[k_p],
\]
where $k_1,k_2,\ldots,k_p \in \bbN$.  

\begin{thm}\mylabel{T:poly}
Given $Q_1,\ldots, Q_p$, the number $\chi_\Phi(Q_{k_1,k_2,\ldots,k_p})$ of totally frustrated states with spins from $Q_{k_1,k_2,\ldots,k_p}$ is given by the multivariate polynomial 
\begin{equation}\mylabel{E:poly}
\chi_{\Phi;Q_1,\ldots,Q_p}(k_1,\ldots,k_p) = 
\sum_{A \in \cH_\Phi}  \mu_{\cH_\Phi}(\eset,A)  \prod_{j=1}^m \bigg[  \sum_{i=1}^p  k_i \,  |\Fix_i \fH_{A_j}(v_j)| \bigg] ,
\end{equation}
where $A_1,\ldots,A_m$ denote the connected components of $A$ and $v_j$ is any vertex of $A_j$.

If not identically $0$, the polynomial has total degree $n$ and the terms of highest degree are those of the expression 
\[
\prod_{v\in V}  \sum_{i=1}^p \  k_i \, \Big[ |Q_i| - \big| \bigcup_{l_v} \Fix_i(\phi(l_v)) \big| \Big] ,
\]
where $l_v$ ranges over all loops incident with $v$.
\end{thm}

Remember that an isolated vertex of $A$ must be treated as a connected component, as its value of $h_Q$ is $\sum_i k_i|Q_i|$.

\begin{proof}
We first give a simple proof of polynomiality, degree, and highest terms, without the explicit formula \eqref{E:poly}.  We use induction on the number of links, employing deletion and contraction.  
Let $\Fix_i(g)$ denote the fixed set of the action of $g$ on $Q_i$.

First consider the case of a single vertex $v$.  For each loop $l_v$, the spins in its fixed set are not allowed to color $v$.  The number of totally frustrated states is therefore 
\begin{align*}
t(v) :&=  \sum_{i=1}^p \  \text{number of spins in $Q_i\times[k_i]$ not fixed by any loop gain at $v$} \\
&=  \sum_{i=1}^p  k_i \Big| Q_i \setm \bigcup_{l_v} \Fix_i(\phi(l_v)) \Big| .
\end{align*}
If $\Phi$ has several vertices, the number of totally frustrated states is the product $\prod_{v\in V} t(v)$, by Proposition \myref{T:mult}.  Thus, the polynomial is homogeneous with total degree $n$, unless there are no totally frustrated states at all.

Now, apply Equation \eqref{E:dc} to a gain graph with a link $e$.  We may assume $\Phi$ does have a totally frustrated state.  We find that $\chi_\Phi(Q)$ is the difference of one polynomial of total degree $n$ and another with total degree $n-1$.  (The former cannot be identically zero, since that would mean $\Phi\setm e$ has no totally frustrated states, hence the same would be true of $\Phi$, contrary to assumption.)  The highest-degree terms of $\chi_\Phi(Q)$, having degree $n$, are those of $\chi_{\Phi\setm e}(Q)$, which by induction are the ones specified in the statement.

The precise formula comes from Theorem \myref{T:detailed}.  It depends on evaluating $h_Q(A)$ for the special spin set $Q$.  Since $\Fix \fH = \bigcup_i \ (\Fix_i \fH)\times[k_i]$, for a connected subgraph with edge set $A_j$ we have
\[
h_Q(A_j) = \sum_{i=1}^p  k_i   |\Fix_i \fH_{A_j}(v_j)| 
\]
where $v_j$ is any vertex of $A_j$.  Thus we immediately obtain \eqref{E:poly}.
\end{proof}

Naturally, the chromatic polynomials are special cases of the grand multivariate polynomial.  
The zero-free chromatic polynomial corresponds to $p=1$ and $Q_1=\fG$ with variable $\lambda:=k_1|Q_1|$, and the chromatic polynomial corresponds to $p=2$, $Q_1=\fG$, $|Q_2|=1$, and $k_2=1$ with variable $\lambda:=k_1|Q_1| + |Q_2|$.

\subsection*{Many zeroes}  The general example with $Q_1=\fG$ and $|Q_2|=1$ is near enough to standard gain-graph coloring to be interesting.  By analogy with standard gain-graph coloring, one might think of $Q_2\times[k_2]$ as the set $\{0\}\times[k_2]$ consisting of $k_2$ different zeroes.  

The fixed sets of the holonomy groups have sizes $h_{Q_2}(A_j) = |Q_2| = 1$ and $$h_{Q_1}(A_j) = |Q_1| \text{ or } 1$$ depending on whether $A_j$ is balanced or unbalanced.  
A connected edge set that is balanced is holonomy closed if and only if it is a maximal balanced set on its vertices.  
A connected edge set that is unbalanced is holonomy closed if and only if it is closed in the underlying graph, i.e., it is a connected induced subgraph.  
Consequently, a set is holonomy closed if and only if each connected component is either an induced subgraph or a maximal balanced set on its vertex set.  These sets include all flats of the frame matroid $G(\Phi)$ and also sets obtained by taking one or more unbalanced components of a flat, partitioning each component's vertex set, and taking the induced subgraphs on the blocks of the partition.  
From Equation \eqref{E:poly}, the formula is
\begin{equation}\mylabel{E:regularfixed}
\begin{aligned}
\chi_{\Phi;Q_1,Q_2}(k_1,k_2) &= 
\sum_{A \in \cH_\Phi}  \mu_{\cH_\Phi}(\eset,A)  \prod_{j=1}^m \big[ k_1 h_{Q_1}(A_j) +  k_2 \big] \\
&= \sum_{A \in \cH_\Phi}  \mu_{\cH_\Phi}(\eset,A) [ k_1 |\fG| +  k_2 ]^{b(A)}  k_2^{c(A)-b(A)} .
\end{aligned}
\end{equation}

\section*{A question}

\subsection*{Negative numbers?}
The striking parallelism with the ordinary theory of graph and gain-graph coloring skips one remarkable feature of the latter theories: the interpretation of the chromatic polynomials at negative arguments $\lambda$.  Can this be repeated for the grand multivariate polynomial?  It is not clear even how to make sense of such a question because there is no one variable that corresponds to $\lambda$.  Yet, it is tantalizing.

%%%%%%%%%%%%%%%%%%%%%%%%


\begin{thebibliography}{9}

\bibitem{BOTutte} Thomas Brylawski and James Oxley, 
The Tutte polynomial and its applications.
In Neil White, ed., \emph{Matroid Applications}, Ch.\ 6, pp.\ 123--225.
Encyc.\ Math.\ Appl., 40.
Cambridge Univ.\ Press, Cambridge, 1992.
MR 93k:05060.  Zbl.\ 769.05026.

\bibitem{DM} Patrick Doreian and Andrej Mrvar, 
A partitioning approach to structural balance.  
\emph{Social Networks} {\bf 18} (1996), 149--168.  

%\bibitem{HPT} David Forge and Thomas Zaslavsky,  Huge polynomial Tutte invariants of weighted gain graphs.  In preparation.
%MR .  Zbl.\ .

\bibitem{FCT} Gian-Carlo Rota,
On the foundations of combinatorial theory:  I.\ Theory of M{\"o}bius functions.
\emph{Z.\ Wahrsch.\ verw.\ Gebiete} {\bf 2} (1964), 340--368.
MR 30 \#4688.  Zbl.\ 121, 24f (e: 121.02406).

\bibitem{RR} S.~S.~Ryshkov and K.~A.~Rybnikov, Jr.,
The theory of quality translations with applications to tilings.
\emph{European J.\ Combin.}\ {\bf 18} (1997), 431--444.
MR 98d:52031. Zbl.\ 881.52015.

\bibitem{CBA} Konstantin Rybnikov and Thomas Zaslavsky,
Criteria for balance in abelian gain graphs, with an application to piecewise-linear geometry.
\emph{Discrete Comput.\ Geom.} {\bf 34} (2005),  no. 2, 251--268.
{\tt arXiv.org} math.CO/0210052.
MR 2006f:05086.  Zbl.\ 1074.05047.

\bibitem{EC1} Richard P.\ Stanley, \emph{Enumerative Combinatorics}, Vol.\ 1.  Wadsworth \& Brooks/Cole, Monterey, Calif., 1986.  MR 87j:05003.  Zbl.\ 608.05001.

%\bibitem{EC2} Richard P.\ Stanley, \emph{Enumerative Combinatorics}, Vol.\ 2.  With an appendix by Sergey Fomin.  Cambridge Stud.\ Adv.\ Math., Vol.\ 62.  Cambridge University Press, Cambridge, 1999.  MR 2000k:05026.  Zbl.\ 928.05001.

\bibitem{WW} Michelle L.\ Wachs and James W.\ Walker, 
On geometric semilattices. 
\emph{Order} {\bf 2} (1986), 367--385.
MR 87f:06004.  Zbl.\ 589.06005.

\bibitem{WelshKnots} D.J.A.~Welsh, 
\emph{Complexity: Knots, Colourings and Counting}.  
London Math.\ Soc.\ Lecture Note Ser., 186.  
Cambridge Univ.\ Press, Cambridge, Eng., 1993.
MR 94m:57027.  Zbl.\ 799.68008.
% No signed edges (except as notational convenience); his edges are in effect all neg.

\bibitem{muChapter} Thomas Zaslavsky, 
The M\"obius function and the characteristic polynomial.  
In: Neil White, ed., \emph{Combinatorial Geometries}, 
Chapter 7, pp.\ 114--138.  
Encyc.\ Math.\  Appl., Vol.\ 29. 
Cambridge Univ.\ Press, Cambridge, 1987.  
MR 88g:05048 (book).  Zbl.\ 632.05017.

\bibitem{BG1} Thomas Zaslavsky,
Biased graphs.  I.\  Bias, balance, and gains.  
\emph{J.\ Combin.\ Theory Ser.\ B} {\bf 47} (1989), 32--52.  
MR 90k:05138.  Zbl.\ 714.05057.

\bibitem{BG2} Thomas Zaslavsky,
Biased graphs.  II.\  The three matroids.  
\emph{J.\ Combin.\ Theory Ser.\ B} {\bf 51} (1991), 46--72.  
MR 91m:05056.  Zbl.\ 763.05096.

\bibitem{BG3} Thomas Zaslavsky,
Biased graphs.  III.\  Chromatic and dichromatic invariants.  
\emph{J.\ Combin.\ Theory Ser.\ B} {\bf 64} (1995), 17--88.
MR 96g:05139.  Zbl.\ 857.05088.

\bibitem{BSG} Thomas Zaslavsky, 
A mathematical bibliography of signed and gain graphs and allied areas.  
\emph{Electronic J.\ Combin.}, Dynamic Surveys in Combinatorics (1998), No.\ DS8 (electronic).
MR 2000m:05001a.  Zbl.\ 898.05001.

\bibitem{SetCol} Thomas Zaslavsky, 
A new distribution problem of balls into urns, and how to color a graph by different-sized sets.
Submitted.

\end{thebibliography}
\end{document}